\documentclass[a4paper,12pt]{article}
\usepackage{mathtext}
\usepackage[T2A]{fontenc}
\usepackage[utf8]{inputenc}
\usepackage[english]{babel}
\usepackage{latexsym,amsfonts,amssymb,amsmath,longtable,amsthm}
\usepackage{cite}
\usepackage[pic]{xy}
\usepackage{bbm,makeidx}
\usepackage[centerlast,small]{caption}
\usepackage{amsmath}
\begin{document}
\begin{center}
{\LARGE The $R_{\infty}$-property for Chevalley groups of types $B_l, C_l, D_l$ over integral domains\footnote{The work is supported by  Russian Science Foundation (project
14-21-00065).
}}

~

~

{\large Timur R. Nasybullov}

~

Sobolev Institute of Mathematics \& Novosibirsk State University

ntr@math.nsc.ru
\end{center}
\newcounter{thelem}
\newtheorem{lem}{{\scshape Lemma}}
\newtheorem{ttt}{{\scshape Theorem}}
\newtheorem{cor}{{\scshape Corollary}}
\newtheorem{prp}[lem]{{\scshape Proposition}}
\renewcommand{\theenumi}{(\asbuk{enumi})}
\renewcommand{\labelenumi}{\asbuk{enumi})}
~~~~~~~~~~~~~~\parbox{11cm}{\begin{center}\textbf{Abstract}
\end{center}
\small We prove that Chevalley groups of the classical series $B_l, C_l, D_l$ over an integral domain of zero characteristic, which has periodic automorphism group, possess the $R_{\infty}$-property.\normalsize}

~

\section{Introduction}

 Let $G$ be a group and $\varphi$ be an automorphism of $G$. Elements  $x, y$ of the group $G$ are said to be  (\emph{twisted}) $\varphi$-\emph{conjugated} ($x\sim_{\varphi}y$) if there exists an element $z\in G$ such  that $x=zy\varphi(z^{-1})$. The relation of $\varphi$-conjugacy is an equivalence relation and it devides the group into $\varphi$-conjugacy classes. The number $R(\varphi)$ of these classes is called the \emph{Reidemeister number} of the automorphism $\varphi$. If $R(\varphi)$ is infinite for any automorphism $\varphi$, then $G$ is said to possess  the $R_{\infty}$-property.

The problem of determining groups which possess the $R_{\infty}$-property was formulated by A. Fel'shtyn and R. Hill \cite{FH}. One of the first general results in this area was obtained by A.~Fel'shtyn, G.~Levitt and M.~Lustig, they proved that non-elementary Gromov hyperbolic groups possess the $R_{\infty}$-property \cite{F,lelu}. Another extensive result was established by A.~Fel'shtyn and E.~Troitsky, they proved that any non-amenable residually finite finitely generated group possesses the $R_{\infty}$-property\cite{FT1}. This wide class of groups contains a lot of finitely generated linear groups, in particular, general linear groups ${\rm GL}_n(\mathbb{Z})$, special linear groups ${\rm SL}_n(\mathbb{Z})$, symplectic groups ${\rm Sp}_{2n}(\mathbb{Z})$. In the paper \cite{Nas2} the author considered some infinitely generated linear groups. In particular, it was proved, that any Chevalley group (of normal type) over an algebraically closed field $F$ of zero characteristic possesses the $R_{\infty}$-property if the transcendence degree of the field $F$ over  $\mathbb{Q}$ is finite.

In this paper we study the $R_{\infty}$-property for Chevalley groups of the classical series $B_l, C_l, D_l$ over integral domains which are not necessarily fields. The main result of the paper is the following
\begin{ttt} Let $G$ be a Chevalley group of type $B_l, C_l$ or $D_l$ over a local integral domain $R$ of zero characteristic. If the automorphism group of the ring $R$ is periodic, then $G$ possesses the $R_{\infty}$ property.
\end{ttt}

In the paper \cite{Nas} similar result was proved for Chevalley groups of type $A_l$, therefore we do not consider the case of  root system $A_l$ in the present paper.

The localization $\mathbb{Z}_{p\mathbb{Z}}$ of the ring of integers $\mathbb{Z}$ by the ideal $p\mathbb{Z}$ is a local integral domain of characteristic zero with the trivial automorphism group and therefore it satisfies the conditions of the theorem.

The condition that the ring $R$ has  characteristic zero is essential. It follows from the result of R.~Steinberg \cite[Theorem 10.1]{S} which says that for any connected linear algebraic group over an algebraically closed field of non-zero characteristic, there always exists an automorphism  $\varphi$ for which $R(\varphi)=1$.

At present, there are no examples of integral domains of characteristic zero  such that Chevalley groups over these domains do not possess the $R_{\infty}$-property. The author believes that it is possible to discard the condition that the automorphism group of the ring $R$ is periodic. The result \cite[Theorem 1]{Nas2} gives a lot of examples of fields of characteristic zero with non-periodic automorphism group such that Chevalley groups over these fields posess the $R_{\infty}$-property.

E.~Jabara studied groups which do not possess the $R_{\infty}$-property. In particular, he proved that any residually finite group which admits an automorphism $\varphi$ of prime order with $R(\varphi)<\infty$ is virtually nilpotent.

\section{Preliminaries}

We use classical notation. Symbols $I_n$ and $O_{n\times m}$ mean the identity $n\times n$ matrix and the $n\times m$ matrix with zero entries, respectively. If  $A$ an $n\times n$ matrix and $B$ an $m\times m$ matrix, then the symbol  $A\oplus B$ denotes the direct sum of the matrices $A$ and $B$, i.~e. the block-diagonal $(m+n)\times (m+n)$ matrix
$$
\newcommand{\tempa}{\multicolumn{1}{c|}{A}}
\newcommand{\tempb}{\multicolumn{1}{|c}{B}}
\begin{pmatrix}
\tempa&O_{n\times m}\\\cline{1-2}
O_{m\times n}&\tempb
\end{pmatrix}
.
$$
  It is obvious that for a pair of $n\times n$ matrices $A_1, A_2$ and for a pair of $m\times m$ matrices $B_1,B_2$ we have
$(A_1\oplus B_1)(A_2\oplus B_2)=A_1A_2\oplus B_1B_2$, $(A_1\oplus B_1)^{-1}=A_1^{-1}\oplus B_1^{-1}$.

 The orthogonal group ${\rm O}_{l}(R,f)$, which preserves a quadratic form $f$, and the symplectic group ${\rm Sp}_{2l}(R)$   over a ring $R$ are defined by the formulas
\begin{eqnarray}
\nonumber &&{\rm Sp}_{2l}(R)=\left\{A\in {\rm GL}_{2l}(R) ~|~A\begin{pmatrix}O_{l\times l}&I_l\\
-I_l&O_{l\times l}
\end{pmatrix}A^T=\begin{pmatrix}O_{l\times l}&I_l\\
-I_l&O_{l\times l}
\end{pmatrix}\right\},\\
\nonumber \\
\nonumber &&{\rm O}_{l}(R,f)=\left\{A\in {\rm GL}_{l}(R) ~|~A[f]A^T=[f]\right\},
\end{eqnarray}
where $[f]$ is a matrix of the quadratic form $f$ and $^T$ denotes transpose. We denote by $\Omega_l(R,f)$ the derived subgroup of ${\rm O}_l(R,f)$. Factoring groups ${\rm Sp}_{2l}(R)$, ${\rm O}_{l}(R,f)$ and $\Omega_l(R,f)$ by their center we obtain the corresponding projective groups  ${\rm PSp}_{2l}(R)$, ${\rm PO}_{l}(R,f)$ and ${\rm P}\Omega_l(R,f)$.

The following proposition about the number of twisted conjugacy classes in a group and in a quotient group was proved in \cite[Lemmas 2.1, 2.2]{MS}.
 \begin{prp}\label{pr3} Let
$$1\rightarrow N\rightarrow G \rightarrow A\rightarrow 1$$
be a short exact sequence of groups, and $N$ be a characteristic subgroup of $G$.
\begin{description}
\item a. ~~If $A$ possesses the $R_{\infty}$-property, then $G$ possesses the $R_{\infty}$-property.
\item b. ~~If $N$ is a finite group and $G$ possesses the $R_{\infty}$-property, then $A$ possesses the $R_{\infty}$-property.
\end{description}
\end{prp}
 The following proposition about the connection between the Reidemeister number of the automorphism $\varphi$ and  the automorphism $\varphi\varphi_H$, where $\varphi_H$ is an inner automorphism induced by the element $H$, can be found in \cite[Corollary 3.2]{feintro} .
\begin{prp}\label{pr1} Let $\varphi, \varphi_H$ be an automorphism and an inner automorphism of the group $G$, respectively. Then $R(\varphi\varphi_H)=R(\varphi)$.
\end{prp}
An associative and commutative ring $R$ is said to be an integral domain if it contains the unit element $1$ and it has no zero devisors. The following simple proposition of ring theory can be found in \cite[Lemma 1]{Nas}
\begin{prp}\label{lem2} Let $K$ be an integral domain and $M$ be an infinite subset of $K$. Then for any polynomial $f$ of non-zero degree the set $P=\{f(a): a \in M\}$ is infinite.
\end{prp}

\section{Proof of the main result}
\textbf{{\scshape Theorem 1}} \emph{Let $G$ be a Chevalley group of  type $B_l, C_l$ or $D_l$ over a local integral domain $R$ of zero characteristic. If the automorphism group of the ring $R$ is periodic, then $G$ possesses the $R_{\infty}$-property.}\\
\textbf{Proof.} We separately consider all the types of root systems.

\emph{Case 1.} \emph{The root system has the type $C_l$.} Since the quotient group $G/Z(G)$ is isomorphic to the elementary Chevalley group $C_l(R)$ \cite[\S12.1]{Car}, then by the proposition \ref{pr3}($a$) it is sufficient to prove that the group $C_l(R)$ possesses the $R_{\infty}$-property.

The group $C_l(R)$ is known to be isomorphic to the projective symplectic group ${\rm PSp}_{2l}(R)$ over the ring $R$  \cite[\S 11.3]{Car}. Since the center of the group ${\rm Sp}_{2l}(R)$ is finite, then by the proposition \ref{pr3}($b$) we can consider $G={\rm Sp}_{2l}(R)$ and prove that this group possesses the $R_{\infty}$-property.

Let $T$ be a variable and $y$ be an element of the ring $R$. Denote by the symbols $X(T)$ and $Y(y)$ the following $2l\times 2l$ matrices
$$ X(T)=\begin{pmatrix}
T\oplus I_{l-1} & I_l \\
-I_l&O_{l\times l}
\end{pmatrix}~~~~~Y(y)=\begin{pmatrix}
I_l & O_{l\times l} \\
O_{l\times l}&yI_l
\end{pmatrix}.$$
Let $Z_y(T)$ be the product of $X(T)$ and $Y(y)$.
 $$Z_y(T)=X(T)Y(y)=\begin{pmatrix}
T\oplus I_{l-1} & yI_l \\
-I_l&O_{l\times l}
\end{pmatrix}$$
By direct calculations we have that for every element $x$ of the ring $R$ the matrix $X(x)$ belongs to $G={\rm Sp}_{2l}(R)$.

Let us prove the following auxiliary statement:\\
For any positive integer $k$ and for every elements $y_1,\dots,y_k$ of the ring $R$ the matrix
$Z_{y_1}(T)\dots Z_{y_k}(T)$ has the form
$$\begin{pmatrix}
f_k(T)\oplus a_kI_{l-1}&g_k(T)\oplus b_kI_{l-1}\\
h_k(T)\oplus c_kI_{l-1}&p_k(T) \oplus d_kI_{l-1}
\end{pmatrix},$$
where $a_k,b_k,c_k,d_k$ are elements of the ring $R$ and $f_k,g_k,h_k,p_k$ are polynomials with coefficients from the ring $R$ such that the degree of $f_k$ is equal to $k$ and degrees of polynomials $g_k,h_k,p_k$ are less than $k$.

  To prove this statement we use induction on the parameter $k$. If $k=1$, then the statement is obvious. Suppose that the statement holds for the number $k-1$, i.~e. the following equality holds
\begin{equation*}
Z_{y_1}(T)\dots Z_{y_{k-1}}(T)=\begin{pmatrix}
f_{k-1}\oplus aI_{l-1}&g_{k-1}\oplus bI_{l-1}\\
h_{k-1}\oplus cI_{l-1}&p_{k-1} \oplus dI_{l-1}
\end{pmatrix},
  \end{equation*}
 where degree of the polynomial $f_{k-1}$ is equal to $k-1$ and degrees of the polynomials $g_{k-1},h_{k-1},p_{k-1}$ are equal to $n,m,r<k-1$, respectively. Then we have
\begin{multline*}
Z_{y_1}(T)\dots Z_{y_k}(T)=\begin{pmatrix}
f_{k-1}\oplus aI_{l-1}&g_{k-1}\oplus bI_{l-1}\\
h_{k-1}\oplus cI_{l-1}&p_{k-1} \oplus dI_{l-1}
\end{pmatrix}
\begin{pmatrix}
T\oplus I_{l-1} & y_kI_l \\
-I_l&O_{l\times l}
\end{pmatrix}
=\\=\begin{pmatrix}
(Tf_{k-1}-g_{k-1}) \oplus (a-b)I_{l-1}&y_kf_{k-1}\oplus y_kaI_{l-1}\\
(Th_{k-1}-p_{k-1})\oplus (c-d)I_{l-1}&y_kh_{k-1}\oplus y_kcI_{l-1}
\end{pmatrix}.
  \end{multline*}
Let us look at degrees of the resulting polynomials. A polynomial in the position $(1,1)$ has the degree $k-1+1=k$; a polynomial in the position $(1,l+1)$ has the degree $k-1<k$; the degree of a polynomial in the  position $(l+1,1)$ is less than or equal to $max\{{\rm deg}(Th(T)),{\rm deg}(p(T))\}=max\{m+1,r\}<max\{k-1+1,k-1\}=k$; and the degree of a polynomial in the position $(l+1,l+1)$ is equal to $m<k-1<k$. The auxiliary statement is proved.
In particular for any positive integer $k$ and for every elements $y_1,\dots,y_k$ of the ring $R$ the trace of the matrix $Z_{y_1}(T)\dots Z_{y_k}(T)$ is a polynomial of degree $k$ with coefficient from the ring $R$.

To prove that the group $G={\rm Sp}_{2l}(R)$ possesses the $R_{\infty}$-property we consider an arbitrary automorphism $\varphi$ of the group $G$ and prove that $R(\varphi)=\infty$. In the papers \cite{MM,Bloh} it is proved that $\varphi$ acts by the rule
$$\varphi: A \mapsto H_1H_2\overline{\delta}(A)H_2^{-1}H_1^{-1},$$
where $\overline{\delta}$ is an automorphism which is induced by the automorphism $\delta$ of the ring $R$
$$\overline{\delta}:A=(a_{ij})\mapsto (\delta(a_{ij})),$$
the matrix $H_1$ belongs to $G$ and the matrix $H_2$ has the form
$$H_2=\begin{pmatrix}
I_l & O_{l\times l} \\
O_{l\times l}&\beta I_l
\end{pmatrix}=Y(\beta)$$
for a certain invertible element $\beta$ of the ring $R$. By the proposition \ref{pr1} we can consider that $\varphi$ acts by the rule
$$\varphi: A \mapsto H_2\overline{\delta}(A)H_2^{-1}.$$

Since an automorphism group of the ring $R$ is periodic, then there exists a number $k$ such that $\overline{\delta}^k=id$. Let $\psi$ be the following function
$$\psi(T)={\rm tr}\left(Z_{\beta}(T)Z_{\delta(\beta)}(T)\dots Z_{\delta^{k-1}(\beta)}(T)\right),$$
 which is a polynomial of the degree $k$ (as we already noted in the auxiliary statement). By the proposition \ref{lem2} there exists an infinite set of integers $a_1,a_2,\dots\in\mathbb{Z}\subseteq R$ such that $\psi(a_i)\neq \psi(a_j) \text{ for } i\neq j$.

Consider the set of matrices $A_1,A_2,\dots$, where $A_i=X(a_i)$, and suppose that $R(\varphi)<\infty$. Then there exist two numbers $i\neq j$ such that $A_i\sim_{\varphi}A_j$, i.~e. for  a certain matrix $D$ the following equality holds
$$A_i=DA_j\varphi(D^{-1})=DA_jH_2\overline{\delta}(D^{-1})H_2^{-1}.$$
If we multiply this equality by the matrix $H_2$ we have
\begin{equation}\label{eq1}
Z_{\beta}(a_i)=A_iH_2=DA_jH_2\overline{\delta}(D^{-1})=DZ_{\beta}(a_j)\overline{\delta}(D^{-1})
\end{equation}
 since $A_iH_2=X(a_i)Y(\beta)=Z_{\beta}(a_i)$.

 Since $\delta$ is an automorphism of the ring $R$,  it acts identically on the subring of integers and therefore $\overline{\delta}(Z_{\beta}(a_i))=Z_{\delta(\beta)}(a_i)$.
Since $\overline{\delta}^k=id$, acting by degrees of the automorphism $\overline{\delta}$ on the equality (\ref{eq1}) we have the following system of equalities
\begin{eqnarray}
\nonumber  Z_{\beta}(a_i) &=& D Z_{\beta}(a_j)\overline{\delta}(D^{-1}), \\
\nonumber  Z_{\delta(\beta)}(a_i) &=& \overline{\delta}(D)Z_{\delta(\beta)}(a_j)\overline{\delta}^{2}(D^{-1}),\\
\nonumber  \vdots &\vdots& \vdots \\
\nonumber  Z_{\delta^{k-1}(\beta)}(a_i) &=& \overline{\delta}^{m-1}(D)Z_{\delta^{k-1}(\beta)}(a_i)D^{-1}. \end{eqnarray}
 If we multiply all of these equalities we conclude that $$Z_{\beta}(a_i)Z_{\delta(\beta)}(a_i)\dots Z_{\delta^{k-1}(\beta)}(a_i)=DZ_{\beta}(a_j)Z_{\delta(\beta)}(a_j)\dots Z_{\delta^{k-1}(\beta)}(a_j)D^{-1},$$
 i.~e. the matrices $Z_{\beta}(a_i)Z_{\delta(\beta)}(a_i)\dots Z_{\delta^{k-1}(\beta)}(a_i)$ and $Z_{\beta}(a_j)Z_{\delta(\beta)}(a_j)\dots Z_{\delta^{k-1}(\beta)}(a_j)$ are conjugated. Therefore, their traces are the same and $\psi(a_i)=\psi(a_j)$. It contradicts to the choice of the elements $a_1,a_2,\dots$ Then the matrices $A_i$ and $A_j$ can not be $\varphi$-conjugated and therefore $R(\varphi)=\infty$.

\emph{Case 2.} \emph{The root system has the type $D_l$.} By the arguments of the case 1 it is sufficient to prove that the elementary Chevalley group $D_l(R)$  possesses the $R_{\infty}$-property.

It is well known that $D_l(R)\cong {\rm P}\Omega_{2l}(R,f_D)$  \cite[\S 11.3]{Car}, where the matrix of the quadratic form $f_D$ has the form
$$[f_D]=\begin{pmatrix}O_{l\times l}& I_l\\
I_l&O_{l\times l}
\end{pmatrix}.$$

Since the center of the group $\Omega_{2l}(R,f_D)$ is finite, then by the proposition \ref{pr3}($b$) we can consider $G=\Omega_{2l}(R,f_D)$ and prove the $R_{\infty}$-property for the group $\Omega_{2l}(R,f_D)$.

Let $T$ be a variable and $X(T),~Y(T)$ be the following matrices
$$X(T)=\begin{pmatrix} 1&T\\
0&1
\end{pmatrix}\oplus I_{l-2}\oplus \begin{pmatrix}1&0\\
-T&1
\end{pmatrix}\oplus I_{l-2},$$
$$Y(T)=\begin{pmatrix} T&1\\
-1&0
\end{pmatrix}\oplus I_{l-2}\oplus \begin{pmatrix}0&1\\
-1&T
\end{pmatrix}\oplus I_{l-2}.$$
Let $Z(T)$ be the commutator  of $X(T)$ and $Y(T)$
$$Z(T)=[X(T),Y(T)]=\begin{pmatrix} T^2+1&-T\\
-T&1
\end{pmatrix}\oplus I_{l-2}\oplus \begin{pmatrix}1&T\\
T&T^2+1
\end{pmatrix}\oplus I_{l-2}.$$
By direct calculations we see that for any element $x$ of the ring $R$ the matrices $X(x), Y(x)$ belong to ${\rm O}_{2l}(R,f_D)$ and therefore $Z(x)$ belongs to $\Omega_{2l}(R,f_D)$.

Let us show that for every positive integer $k$ the trace of the matrix $Z(T)^k$ is a non-constants integral polynomial. To do it we prove more general result: For any positive integer $k$ the matrix $Z(T)^k$ has the form
$$\begin{pmatrix} f_k(T)&g_k(T)\\
h_k(T)&p_k(T)
\end{pmatrix}\oplus I_{l-2}\oplus \begin{pmatrix}p_k(T)&-h_k(T)\\
-g_k(T)&f_k(T)
\end{pmatrix}\oplus I_{l-2},$$
where $f_k$ is a polynomial of degree $2k$, and $g_k,h_k,p_k$ a polynomials of degrees which are less than $2k$.

We use induction on the parameter $k$. The basis of induction ($k=1$) is obvious.
Suppose that this statement holds for the number $k-1$, i.~e. for certain integral polynomials $f_{k-1},g_{k-1},h_{k-1},p_{k-1}$ the following equality holds
$$Z(T)^{k-1}=\begin{pmatrix} f_{k-1}&g_{k-1}\\
h_{k-1}&h_{k-1}
\end{pmatrix}\oplus I_{l-2}\oplus \begin{pmatrix}p_{k-1}&-h_{k-1}\\
-g_{k-1}&f_{k-1}
\end{pmatrix}\oplus I_{l-2},$$
 where degree of the polynomial $f_{k-1}$ is equal to $2(k-1)$ and degrees of the polynomials $g_{k-1}$, $h_{k-1}$, $p_{k-1}$ are equal to $n,m,r<2(k-1)$, respectively. Then the matrix $Z(T)^k=Z(T)^{k-1}Z(T)$ has the form
\begin{multline*}
\begin{pmatrix} (T^2+1)f_{k-1}-Tg_{k-1}&-Tf_{k-1}+g_{k-1}\\
(T^2+1)h_{k-1}-Tp_{k-1}&-Th_{k-1}+p_{k-1}
\end{pmatrix}\oplus I_{l-2}\oplus\\
\oplus \begin{pmatrix}-Th_{k-1}+p_{k-1}&-(T^2+1)h_{k-1}+Tp_{k-1}\\
Tf_{k-1}-g_{k-1}&(T^2+1)f_{k-1}-Tg_{k-1}
\end{pmatrix}\oplus I_{l-2}
\end{multline*}
A polynomial in the position $(1,1)$ of this matrix has the degree $2(k-1)+2=2k$; a polynomial in the position $(1,2)$ has the degree $2(k-1)+1=2k-1<2k$; the degree of a polynomial in the position $(2,1)$ is less than or equal to
$$max\{deg((T^2+1)h_{k-1}), deg(Tp_{k-1})\}=max\{m+2,r+1\}<max\{2k,2k-1\}=2k;$$
and the degree of a polynomial in the position $(2,2)$ is less than or equal to
$$max\{deg(Th_{k-1},deg(p_{k-1}))\}=max\{m+1,r\}<max\{2k-1,2k-2\}=2k-1.$$
The auxiliary statement is proved. As a corollary we have that for every positive integer $k$ the function $\psi_k(T)=tr(Z(T)^k)$ is a non-constant integral polynomial.

To prove that the group $G=\Omega_{2l}(R,f_D)$ possesses the $R_{\infty}$-property we consider an arbitrary automorphism $\varphi$ of the group $G$ and prove that $R(\varphi)=\infty$. In the papers \cite{Bun1,Bun2} it is proved that there exist
\begin{description}
\item 1. An inner automorphism $\varphi_H$
$$\varphi_H:A\mapsto HAH^{-1}$$
\item 2. A central automorphism $\Gamma$
$$\Gamma:A\mapsto \gamma(A)A,$$
where $\gamma$ is a homomorphism from the group $G$ into its center $Z(G)$.
\item 3. A ring automorphism $\overline{\delta}$
$$\overline{\delta}:A=(a_{ij})\mapsto (\delta(a_{ij})),$$
where $\delta$ is an automorphism of the ring $R$
\end{description}
such that $\varphi=\varphi_H\Gamma\overline{\delta}$. By the proposition \ref{pr1} we can consider $\varphi=\Gamma\overline{\delta}$.

Since an automorphism group of the ring $R$ is periodic, there exists such a number $k$ that $\overline{\delta}^k=id$. By the proposition \ref{lem2} there exists an infinite set of elements $a_1,a_2,\dots \in \mathbb{Z}\subseteq R$ such that $(\psi_k(a_i))^2\neq(\psi_k(a_j))^2$ for $i\neq j$.

 Consider the set of matrices $A_1=Z(a_1), A_2=Z(a_2),\dots $ and suppose that $R(\varphi)<\infty$. Then there are two $\varphi$-conjugated matrices in the set $A_1, A_2,\dots$, i.~e. for some indexes $i\neq j$ and for some matrix $D\in G$ the following equality holds
\begin{equation}\label{eq2}
A_i=DA_j\varphi(D^{-1})=DA_j\Gamma \overline{\delta}(D^{-1})=DA_jC_1 \overline{\delta}(D^{-1}),
\end{equation}
where the matrix $C_1$ belongs to $Z(G)$.

Since the matrices $A_i, A_j$ have integer coefficients and the automorphism $\delta$ acts identically on the subring of integers, then  $\overline{\delta}(A_i)=A_i$, $\overline{\delta}(A_j)=A_j$. Acting by degrees of the automorphism $\overline{\delta}$ on the equality  (\ref{eq2}) we have the following system of equalities:
\begin{eqnarray}
\nonumber  A_i &=& D A_jC_1\overline{\delta}(D^{-1}), \\
\nonumber  A_i &=& \overline{\delta}(D)A_jC_2\overline{\delta}^{2}(D^{-1}),\\
\nonumber  \vdots &\vdots& \vdots \\
\nonumber  A_i &=& \overline{\delta}^{k-1}(D)A_jC_kD^{-1}.
\end{eqnarray}
If we multiply all of this equalities denoting $C=C_1C_2\dots C_k$, then we have
$$
A_i^k=DCA_j^kD^{-1},
$$
i.~e. the matrices $A_i^k$ and $CA_j^k$ are conjugated and therefore $tr(A_i^k)=tr(CA_j^k)$.
 Since $C\in Z(\Omega_{2l}(R,f_D))=\{ \pm I_{2l}\}$, we have
 $$\psi_k(a_i)=tr(A_i^k)=\pm tr(A_j^k)=\pm\psi_k(a_j).$$
 It contradicts to the choice of the elements $a_1,a_2,\dots$

\emph{Case 3.} \emph{The root system has the type $B_l$.} The elementary Chevalley group $B_l(R)$ is isomorphic to the group ${\rm P}\Omega_{2l+1}(R,f_B)$ \cite[\S 11.3]{Car}, where the matrix of the quadratic forms $f_B$ has the following form
$$[f_B]=1\oplus\begin{pmatrix}O_{l\times l}& I_l\\
I_l&O_{l\times l}
\end{pmatrix}.$$
Using this fact, the proof of the case 3 literally repeats the proof of the case 2 after changing the matrix $Z(T)$ by the matrix $1\oplus Z(T)$, and using the result \cite{Bun3} (instead of \cite{Bun1,Bun2}) about the automorphism group of the Chevalley groups of the type $B_l$.
Theorem is proved.\\[2cm]

\end{document}